\documentclass[review]{elsarticle}

\usepackage{amsmath, amsfonts}
\usepackage{fullpage}
\usepackage{float}
\usepackage{graphicx}
\usepackage{subcaption}
\usepackage{enumitem}
\usepackage{epstopdf}
\usepackage{stackengine}
\usepackage[font=footnotesize, labelfont=bf]{caption}

\biboptions{sort&compress}

\let\temp\varphi%
\let\varphi\phi%
\let\phi\temp%

\newcommand{\sn}{\text{sn}}
\newcommand{\bfK}{\mathbf{K}}
\newcommand{\bfF}{\mathbf{F}}
\newcommand{\bfE}{\mathbf{E}}
\newcommand{\bfPi}{\mathbf{\Pi}}

 \newcommand{\substit}[3]{\left. {#1}%
_{\stackunder[5pt]{}{}}%
 \right|_{%
 \shortstack[l]{{\scriptsize$#2$}\\ {\vspace*{-5pt}\scriptsize$#3$} } }}

\journal{Journal of Sound and Vibration}

\begin{document}

\begin{frontmatter}

\title{Approximation of potential function in the problem of forced escape}

\author{Pavel Kravetc\corref{mycorrespondingauthor}}
\cortext[mycorrespondingauthor]{Corresponding author}
\ead{pmkravets@gmail.com}

\author{Oleg Gendelman}
\ead{ovgend@technion.ac.il}
\address{Faculty of Mechanical Engineering, Technion --- Israel Institute of Technology, 3200003 Haifa, Israel}

\begin{abstract}
\noindent The paper addresses an escape of a classical particle from a potential well under harmonic forcing.  Most dangerous/efficient escape dynamics reveals itself in conditions of 1:1 resonance and can be described in the framework of isolated resonant (IR) approximation. The latter requires reformulation of the problem in terms of action-angle (AA) variables, available only for a handful of the model potentials. The paper suggests approximation of realistic generic potentials by low-order polynomial functions, admissible for the AA transformation, with possible truncation. To illustrate the idea, we first formulate the AA transformation and solve the escape problem in the IR approximation for a generic quartic potential. Then, the model problem for dynamic pull-in in microelectromechanical system (MEMS) is analyzed. The model electrostatic potential is approximated by the quartic polynomials (globally and locally), and quality of predicting the escape thresholds is assessed numerically. Most accurate predictions are delivered by global $L^2$-optimal heuristic approximation.
\end{abstract}


\end{frontmatter}


\section*{Introduction}
\noindent The problem of escape from a potential well under the influence of external forcing, or simply the escape problem, is a well-known problem in both science and engineering. It is often employed to describe transient processes and phenomena such as gravitational collapse, energy harvesting~\cite{mann2009energy}, particle absorption, physics of Josephson junctions~\cite{barone1982physics}, resonance dynamics of oscillatory systems~\cite{quinn1997}, dynamic pull-in in microelectromechanical systems (MEMS)~\cite{alsaleem2010experimental}, and even ship capsizing~\cite{belenky2007,  thompson1992mechanics, virgin1989}, to mention a few. The escape problem dates back to a seminal work by Kramers on thermal activation of chemical reactions~\cite{kramers1940}, where he considered escape under the action of Brownian motion. Even after more than 70 years of active development, this research field remains active nowadays, and contains many open problems~\cite{talkner2012new}.

One encounters the opposite limiting case, if the forcing contains only one Fourier component. In this case, the most salient phenomenon is a resonant escape under the influence of harmonic external force.
Most of the approaches to the problem rely on the numerical methods. However, in recent years an analytic technique --- approximation of the isolated resonance --- was proposed in~\cite{gendelman2018escape}. This method treats the principal $1:1$ resonance though the canonical action-angle (AA) transformation followed by the averaging over the fast phases. As a result one obtains a slow evolution equations of averaged action possessing a first integral which defines a family of resonance manifolds (RM). Initial conditions define a special phase trajectory on the RM. For the zero initial conditions, such special trajectories are called limiting phase trajectories (LPT). There are two mechanisms of escape: saddle mechanism and maximum mechanism. The former corresponds to a passage of the LPT (or other phase trajectory for nonzero IC) through a saddle on the RM, while the latter corresponds to the LPT tangentially crossing the escape barrier. The competition between the two mechanisms yields a theoretical prediction for the critical forcing needed for the escape at a given frequency. The obtained curve features a dip shape with a sharp minimum at the resonance frequency. This method has been proven effective for a variety of potentials including an infinite range potential~\cite{gendelman2018escape} and particular cases of polynomial potentials~\cite{gendelman2019basic, Farid2021}. 

The AA transformation can be performed rigorously only for a handful of the model potentials. To overcome this restriction. in the present work we study the escape from a potential well described by a general quartic polynomial with a two-fold purpose. First of all,  it is the most general case of polynomial potentials for which transformation to AA variables can be done in terms of well-known elliptic functions. Then, we conjecture that forth order polynomial can serve as a good approximation for more intricate potential functions such as, for example, electrostatic potential. To prove the approximation useful we apply it to the escape dynamics in a simple MEMS device --- parallel-plate electrostatic actuator. The escape with or without external forcing is an intrinsic feature systems which combine electrostatic and mechanical forces. In the context of MEMS the escape is a structural instability called pull-in~\cite{younis2011mems}. In particular, a pull-in occurring under the influence of external forcing (e.g, AC loading) is called {\it dynamic}. A potential describing a MEMS actuator contains a singularity which corresponds to the collapse of the plates. Unfortunately, analytical treatment of the transient escape dynamics in such potentials poses a difficult if not an impossible challenge, therefore, finding an appropriate approximation is a great interest to engineers. In this work, we discuss two different approaches to the problem: a global and a local approximation. The former is an ad hoc approach to approximate a given potential with a help of a handful parameters or by fitting the forth-order curve. The latter corresponds to an approximation using a Taylor's polynomial near the minimum of the potential. 

The paper is organized as follows. In Section~\ref{sec:method} we briefly formulate the general problem and outline the method. In Section~\ref{sec:model:quartic} we apply the method to the model with quartic potential, the essential building block for the approximation of more intricate potentials. In Section~\ref{sec:mems} we test different approximation techniques to the model electrostatic potential and discuss their applicability and drawbacks. Finally, Appendix contains the derivations of the main formulae.

\section{Problem formulation and an outline of the method}\label{sec:method}
\noindent The analytic approach has been proposed in~\cite{gendelman2018escape} and used in some  subsequent publications~\cite{gendelman2019basic, Farid2021}. For the convenience of the reader we outline the method here.

\subsection{Formulation of the problem}
\noindent
Let $q$ denote the displacement of a SDOF classical particle of the unit mass which is placed at a local minimum $q=q_0$ of a potential $V(q)$ and is subject to an external harmonic forcing with amplitude $F$, frequency $\Omega$ and phase $\psi$. Without loss of generality one can assume $q_0=0$. Then, the equation of motion of the particle is
\begin{equation}
\label{eq:main}
\ddot{q}(t) + \frac{\mathrm{d} V}{\mathrm{d} q} = F\sin{\left(\Omega t + \psi\right)}.
\end{equation}
A common definition of escape is
$$\lim_{t\to\infty} q(t) \not\in (q_{\text{low}}, q_{\text{high}})$$
where $q_{\text{low}}$ and $q_{\text{high}}$ are lower and upper boundaries of the potential well, respectively.
However, this definition is problematic to use in context of considered problem. First of all, it is impossible to utilize it as an escape criterion in the numerical simulations. Secondly, in some cases (e.g., double-well potential) the aforementioned definition is inapplicable altogether, as according to the Poincar\'e Recurrence Theorem the particle will visit any arbitrary set infinitely many times, and hence, the escape will never happen.
Therefore, we adopt the ``first-hitting" definition instead, i.e., we say that escape occurs if either
$$
\min_{t}{\left\{q(t)\right\}} < q_{\text{low}} \qquad\text{or}\qquad \max_{t}\left\{q(t)\right\} > q_{\text{high}}.
$$
The thresholds $q_{\text{low}}$, $q_{\text{high}}$ can be defined via the maximum energy level $E_{\text{thres}}$ as the solutions to the equation $V(q) = E_{\text{thres}}$.

Alternatively, one can utilize the so-called energy criterion:
$$\max_{t}\left\{E(t)\right\} > E_{\text{thres}},$$
where $E(t)={\dot{q}(t)}^2/2+V(q(t))$ is the total energy of the system.

The central question to the forced escape problem can be formulated in the following way: {\em for a given frequency $\Omega$ in the vicinity of the primary resonance what is the minimal amplitude $F_{\text{crit}}$ needed to trigger an escape?}

\subsection{Method}
\noindent
Equation~\eqref{eq:main} can be rewritten in the Hamiltonian form:
\begin{equation}
\dot q = \frac{\partial H}{\partial p},\qquad \dot p = -\frac{\partial H}{\partial q}
\end{equation}
where the Hamiltonian
\begin{equation}
\label{eq:hamiltonian}
H(q, p) = H_0(p, q) - q F \sin{\left(\Omega t + \psi\right)}
\end{equation} and
\begin{equation}
\label{eq:hamiltonian0}
H_0(p, q) = \frac{p^2}{2} + V(q).
\end{equation}
The basic  Hamiltonian $H_0$ describes the free motion of the particle in the potential well~$V(q)$. We perform a canonical action-angle (AA) transformation using well-known formulae~\cite{Landau1976Mechanics}:
\begin{equation}\label{eq:aa}
I = \frac{1}{2\pi}\oint\limits_{\Gamma_E} p(q,\, E)\mathrm{d}q, \qquad \theta = \frac{\partial}{\partial I} \int\limits_0^q p(x, I) \mathrm{d}x
\end{equation} where $\Gamma_E$ is a phase curve defined by a level set $\left\{H_0=E\right\}$. The canonical transformation does not depend on time explicitly, therefore, the Hamiltonian~\eqref{eq:hamiltonian} can be rewritten in the AA variables:
\begin{equation}
\label{eq:hamil-aa}
H(I, \theta) = H_0(I) - q(I, \theta) F \sin(\Omega t + \psi).
\end{equation} Due to the $2\pi$-periodicity of the angle variable $\theta$, it can be expanded in terms of Fourier series:
\begin{gather}
H = H_0(I) + \frac{\mathrm{i} F}{2} \sum\limits_{m=-\infty}^{\infty} q_m(I) \left(e^{\mathrm{i}(m\theta + \Omega t + \psi)}-e^{\mathrm{i}(m\theta - \Omega t - \psi)}\right),\\
q_m(I) = \bar{q}_{-m}(I).
\end{gather}
Here, $\bar{q}$ denotes the complex conjugation of $q$. The corresponding Hamilton equations are
\begin{eqnarray*}
\dot{I} = -\frac{\partial H}{\partial \theta} = \frac{F}{2}\sum_{m=-\infty}^{\infty} m\, q_m(I)\left(e^{\mathrm{i}(m\theta + \Omega t + \psi)}-e^{\mathrm{i}(m\theta - \Omega t - \psi)}\right),\\
\dot{\theta} = \frac{\partial H}{\partial I} = \frac{\mathrm{d} H_0}{\mathrm{d} I} + \frac{\mathrm{i} F}{2}\sum_{m=-\infty}^{\infty}\frac{\mathrm{d} q_m}{\mathrm{d} I} \left(e^{\mathrm{i}(m\theta + \Omega t + \psi)}-e^{\mathrm{i}(m\theta - \Omega t - \psi)}\right).
\end{eqnarray*}

We consider the primary $1:1$ resonance, i.e., we select $\vartheta = \theta-\Omega t-\psi$ to be the slow phase, and assume all other combinations to be fast. After averaging over the fast phases, we arrive at the following system of slow-flow equations:
\begin{equation}\label{eq:slow:flow}
\begin{aligned}
\dot{J} =& -\frac{F}{2}\left(q_1(J)e^{\mathrm{i}\vartheta}+\bar{q}_1(J)e^{-\mathrm{i}\vartheta}\right),\\
\dot{\vartheta} =& \frac{\partial H_0(J)}{\partial J} - \frac{\mathrm{i} F}{2}\left(\frac{\partial q_1(J)}{\partial J}e^{\mathrm{i} \vartheta}+\frac{\partial \bar{q}_1(J)}{\partial J}e^{-\mathrm{i} \vartheta}  \right) - \Omega,
\end{aligned}
\end{equation} where $J = \langle I(t)\rangle$ denotes the average of the action variable over the fast phases. It is easy to check by differentiation that system~\eqref{eq:slow:flow} possesses the following conservation law:
\begin{equation}
\label{eq:cons:law}
H_0(J) - \frac{\mathrm{i} F}{2} \left(q_1(J) e^{\mathrm{i}\vartheta}-\bar{q}_1(J) e^{-\mathrm{i}\vartheta}\right) - \Omega J = \text{const}.
\end{equation} Often, it is impossible to obtain expression~\eqref{eq:cons:law} in a closed form. However, in order to analyze the escape dynamics, it is sufficient to parameterize~\eqref{eq:cons:law} using averaged energy $\xi = \langle E(t)\rangle$ instead of the averaged action~$J$. In this case, the first integral is
\begin{equation}
\label{eq:cons:law:xi}
C(\vartheta,\,\xi) = \xi - \frac{\mathrm{i} F}{2} \left(q_1(\xi) e^{\mathrm{i}\vartheta}-\bar{q}_1(\xi) e^{-\mathrm{i}\vartheta}\right) - \Omega J(\xi) = C_0.
\end{equation}

Equation \eqref{eq:cons:law:xi} defines a family of $1:1$ resonance manifolds (RMs) on the phase cylinder $(\vartheta,\xi)$. Constant $C_0$ is defined by the initial conditions on the RM, i.e., the values of averaged action $J$ and the slow phase $\vartheta$ at which the system is captured by the RM. We are interested in the escape from the zero initial conditions, hence, the corresponding constant $C_0=0$. This trajectory is often called a {\em limiting phase trajectory} (LPT). Therefore, escape of the particle from the potential well occurs if the LPT reaches the circle~$\xi=E_\text{thres}$.

Based on the behavior of the LPT with varying amplitude $F$ of the external forcing, there are two distinct mechanisms of transition to the escape. The first mechanism is called {\em maximum mechanism} (MM) and it works as follows. At $F = F_{\text{crit}}$, the LPT is tangent to the circle $\xi = E_{\text{max}}$ at some $\vartheta = \vartheta^*$. For $F < F_\text{crit}$, the LPT does not reach the circle $\xi = E_{\text{max}}$. For $F > F_\text{crit}$, the LPT reaches the circle $\xi = E_{\text{max}}$ and the escape occurs. In order to find $F_{\text{crit}}^{\text{MM}}(\Omega)$, i.e.,critical force $F_\text{crit}$ at a given frequency value $\Omega$ for the maximum mechanism, one can solve the equation
\begin{equation}
C(\vartheta^*,\, E_{\text{max}}) = 0,
\end{equation}
where $\vartheta^*$ is defined by equation
 \begin{equation}
\substit{\frac{\partial C}{\partial \vartheta}} {\vartheta=\vartheta^*}{\xi=E_{\text{max}} } =0
\end{equation}

Another way to escape is called the saddle mechanism. It corresponds to the scenario where at $F = F_\text{crit}$, the LPT passes through a saddle point $S = \left(\xi^\dagger, \vartheta^\dagger\right)$, at $F < F_\text{crit}$ the LPT is below point $S$ keeping the particle in the well, and at $F > F_{\text{crit}}$ the LPT connects the circles $\xi = 0$ and $\xi = E_{\text{max}}$, thus, producing the escape trajectory. The saddle point $S$ is defined by the following system
\begin{equation}\label{eq:saddle:condition}
\substit{\frac{\partial C}{\partial \vartheta}}{\vartheta=\vartheta^\dagger}{\xi=\xi^\dagger} = 0,\quad \substit{\frac{\partial C}{\partial \xi}}{\vartheta=\vartheta^\dagger}{\xi=\xi^\dagger} = 0, \quad \det{\left(\mathbf{J}(\nabla C(\vartheta^\dagger, \xi^\dagger))\right)} < 0,
\end{equation} where $\nabla C$ denotes the gradient of $C$ and $\mathbf{J}$ is the Jacobian matrix. The first equation in~\eqref{eq:saddle:condition} immediately yields possible values of $\vartheta^\dagger$. However, the second equation in~\eqref{eq:saddle:condition} is usually a very cumbersome expression. One can avoid dealing with it entirely by using equation of the LPT
\begin{equation}
C(\vartheta, \xi) = 0.
\end{equation} Thus, by solving linear system
\begin{equation}\label{eq:saddle:general}
\substit{\frac{\partial C}{\partial \xi}}{\vartheta= \vartheta^\dagger}{\xi=\xi^\dagger} = 0,\quad C\left(\vartheta^\dagger,\, \xi^\dagger\right) = 0,
\end{equation}
one can obtain a curve $F_{\text{crit}}^{\text{SM}}$ in the space $(\Omega, F)$, parameterized by $\xi^\dagger$. A part of this curve corresponds to the critical forcing needed for the escape at the given frequency~$\Omega$.

The line $F_{\text{crit}}^{\text{MM}}$ and the curve $F_{\text{crit}}^{\text{SM}}$ intersect at a sharp minimum forming a familiar dip shape.

\section{Model with quartic potential}\label{sec:model:quartic}
\noindent A particular case of quartic potential without cubic terms (quadratic-quartic function) was already considered in~\cite{gendelman2019basic}. Due to an additional symmetry of this potential, the conservation law~\eqref{eq:cons:law} as well as the AA transformation together with its inverse, can be elegantly expressed in closed forms. In this paper, we apply the method for a general quartic potential function. We distinguish two cases of quartic potential: {\it double well potential} and {\it inverted quartic potential}. In other words, we consider potential
\begin{equation}\label{eq:potential}
    V(q)=\frac{1}{2}q^2 + \frac{\alpha}{3} q^3 + \frac{\beta}{4}q^4,
\end{equation} where parameters $\alpha$, $\beta$ are such that, $V(q)$ is
\begin{enumerate}[label={Case \Roman*:}]
\item a double-well potential, i.e. the parameters $\alpha$, $\beta$ satisfy two simple inequalities
\begin{equation}\label{eq:cond:case:1}
\alpha < 0 \quad \text{and}\quad 0 < \beta < \frac{2\alpha^2}{9},
\end{equation}
\item an inverted quartic potential, in which case $\alpha$ can be any real number and $\beta < 0$.
\end{enumerate}

\begin{figure}[H]
\centering
\begin{subfigure}[t]{0.32\textwidth}
\includegraphics*[width=0.9\textwidth]{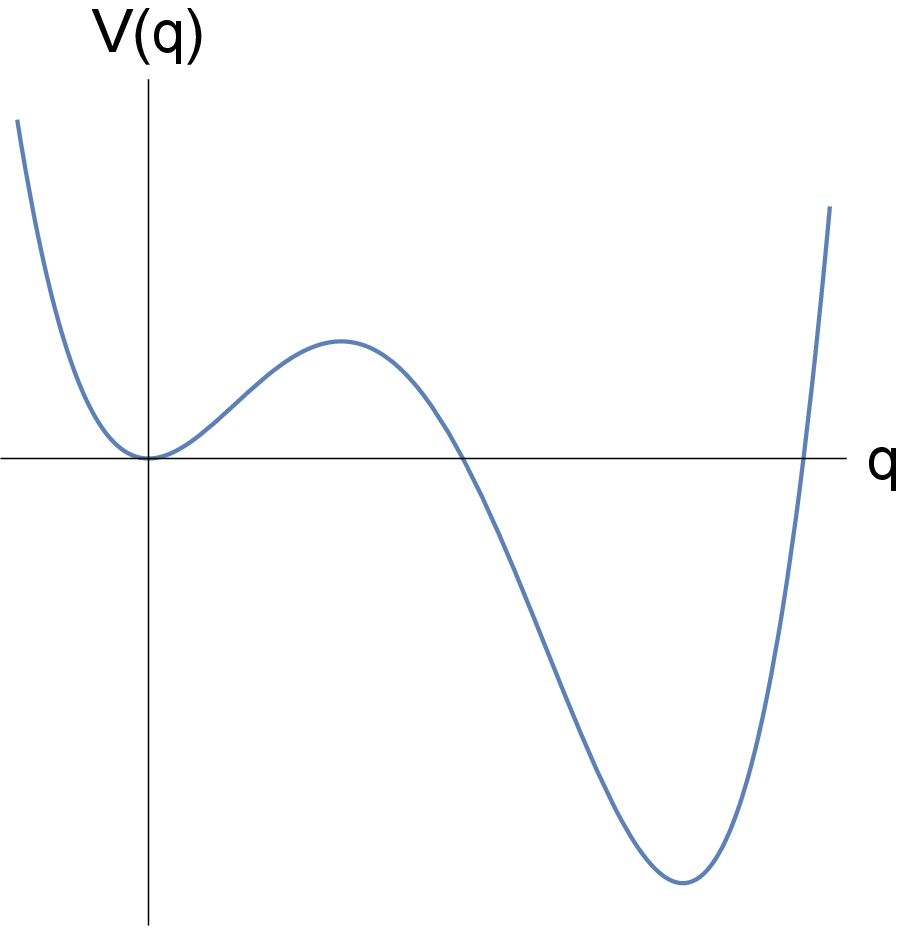}
\caption{\label{fig:case1:pot}}
\end{subfigure}
\hspace{30pt}
\begin{subfigure}[t]{0.32\textwidth}
\includegraphics*[width=0.9\textwidth]{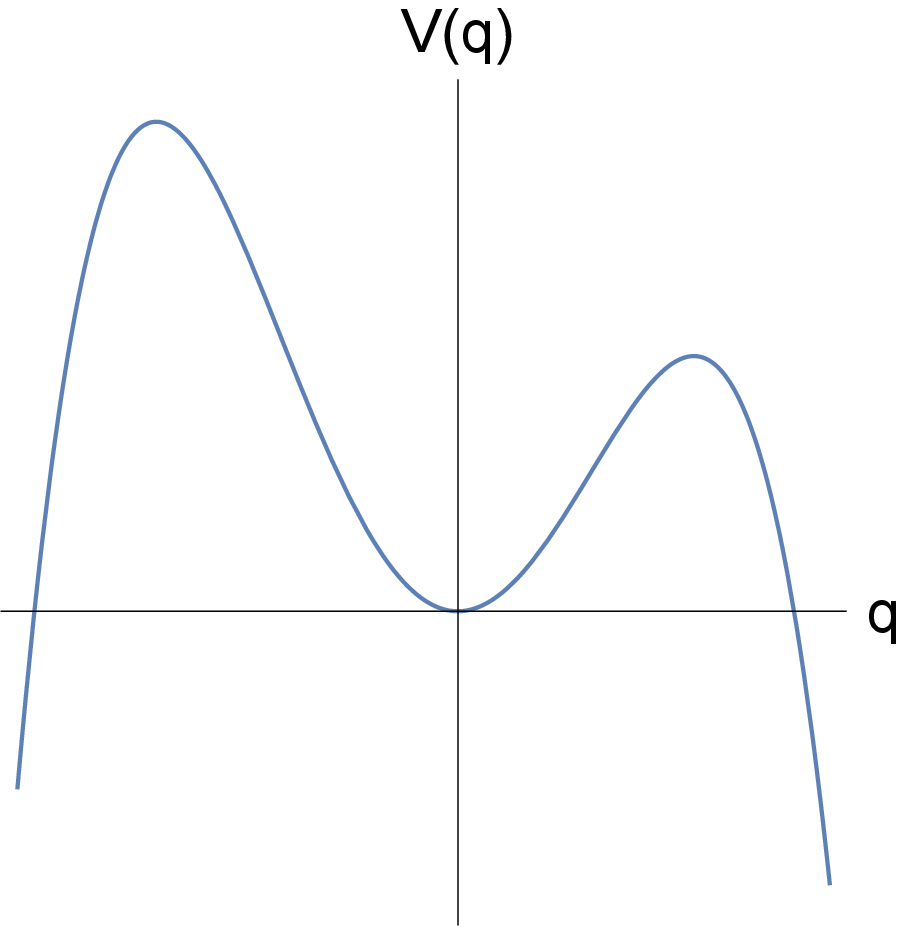}
\caption{\label{fig:case2:pot}}
\end{subfigure}
\caption{Typical shapes of the quartic potential~\eqref{eq:potential}. Panel (a): Case I, $\alpha=-1/2$, $\beta=1/20$; panel (b): Case II, $\alpha=-3/50$, $\beta=17/250$.\label{fig:cases}}
\end{figure}

Case I describes an escape from the shallow well into the deep one, i.e., a transition from a metastable state to the state of the least energy. Escape in the opposite direction is out of scope of the present work, as in this case the derivation of equation~\eqref{eq:cons:law:xi} is too cumbersome. Typical examples of both cases are illustrated on Figure~\ref{fig:cases}.

The difference between two considered cases is rather technical. In particular, different limits of integration in~\eqref{eq:aa} yield slightly different expressions for the conservation law~\eqref{eq:cons:law}. However, their structures are virtually the same.

\subsection{Case I: Double-well Potential}
\noindent Given parameters $\alpha$, $\beta$ satisfy condition \eqref{eq:cond:case:1} the potential $V(q) \to \infty$ as $q \to \pm \infty$, and it has two minima at $q=0$ and $q=\frac{-\alpha + \sqrt{\alpha^2-4\beta}}{2\beta}$ as well as a maximum at 
\begin{equation}\label{eq:q:thres}
q_{\text{thres}}=\frac{-\alpha - \sqrt{\alpha^2-4\beta}}{2\beta}.
\end{equation} If we consider a full (non-truncated) potential well, then the threshold energy level corresponds to the local maximum
\begin{equation}
\label{eq:Emax}
    E_{\text{max}} = V(q_{\text{thres}}) = \frac{\left(\sqrt{\alpha ^2-4 \beta }+\alpha \right)^2 \left(6 \beta-\alpha  \left(\sqrt{\alpha^2-4 \beta }+\alpha \right)\right)}{96 \beta ^3}.
\end{equation}
However, one can choose to use a truncated potential, i.e., to select a cutoff energy level $E_{\text{thres}} < E_{\text{max}}$. Regardless, in case $V(q)$ is a double well potential, the conservation law~\eqref{eq:cons:law:xi} is
\begin{equation}\label{eq:cons:law:case1}
C(\vartheta, \xi) = \xi -FG\sin{\vartheta} - \Omega J = \text{const},
\end{equation} where $\xi=\langle E(t)\rangle$ is the averaged energy and function $G$ is
\[
G(\xi) = \frac{\pi\sqrt{(b-d)(a-c)} \sinh(2\omega)}{2\bfK(k)\sinh(2\omega_0)}.
\] Functions $a = a(\xi)$, $b = b(\xi)$, $c = c(\xi)$, $d = d(\xi)$ are roots of the following forth-order polynomial
\[
\xi - V(q) = \xi - \frac{1}{2} q^2 - \frac{\alpha}{3}q^3 - \frac{\beta}{4}q^4,
\] in descending order ($a>b>c>d$), and
\[
\omega = \frac{\pi \left(\bfK\left(k^\prime\right)-\text{cn}^{-1}\left(\sqrt{\frac{c-d}{a-d}},\, k^\prime\right)\right)}{2 \bfK(k)},\qquad \omega_0 = \frac{\pi \bfK\left(k^\prime\right)}{2 \bfK(k)}, \qquad k^\prime = \sqrt{1-k^2}.
\] The averaged action $J=J(\xi)$ is
\begin{align*}
J = &\frac{1}{48\pi}\sqrt{\frac{2\beta}{(a-c)(b-d)}}\left[(a-c)(b-d)\left(\frac{16 \left(\alpha ^2-3 \beta \right)}{3 \beta ^2}\right)\mathbf{E}(k) +\right.\\ 
&(a-c) (a-d) \left(3 a^2-6 a b-b^2+4 b (c+d)-3 c^2+2 c d-3 d^2\right)\mathbf{K}(k) + 3(a-d)\times\\&\left. \left(-3 a^3+\frac{16 \alpha ^2 a}{9 \beta ^2}-\frac{4 a (a \alpha +3)}{3 \beta}-b^3+(c+d) \left(b^2-(c-d)^2\right)+b \left(c^2+d^2\right)\right)\bfPi(\gamma^2, k)\right].
\end{align*}It is easy to see that $a, b, c, d \in \mathbb{R}$ are well-defined for $0<\xi<E_{\text{max}}$. Although, it is possible to find these roots in exact form using the well-known Ferrari method, the resulting expressions are too awkward to handle. Functions $\mathbf{K}(k)$, $\mathbf{E}(k)$, $\bfPi(\gamma^2,\, k)$ are complete elliptic integrals of the first, the second and the third kind, respectively, with modulus 
\[
k = \sqrt{\frac{(a-b) (c-d)}{(a-c) (b-d)}},
\] and parameter
\[
\gamma^2 = \frac{d-c}{a-c} < 0.
\]
For the derivation of the conservation law~\eqref{eq:cons:law:case1} see Appendix. Equation~\eqref{eq:cons:law:case1} defines the family of the RMs on the phase cylinder~$(\vartheta,\xi)$. Recall that escape from the zero initial conditions occurs when the LPT reaches the circle $\xi = E_{\text{thres}}$.

Equation 
\[\frac{\partial C}  {\partial \vartheta} = 0\]
yields two solutions,
\[
\vartheta^{\dagger}=\frac{\pi}{2} \quad\text{and}\quad \vartheta^{*}=\frac{3\pi}{2}.
\]

By a simple topological argument one can easily show that one of the obtained critical points is a saddle.

Equations~\eqref{eq:saddle:general} then become
\begin{align*}
F G^{\prime}+\Omega J^{\prime}&= 1, \\
F G + \Omega J &=\xi,
\end{align*} thus, expressions
\begin{equation}\label{eq:saddle:case:one}
F(\xi) = \frac{J - \xi J^\prime}{JG^\prime-G J^\prime}, \qquad \Omega(\xi) = \frac{\xi\left(G^\prime- J^{\prime}\right)}{JG^\prime-G J^\prime}
\end{equation}
defines a parametric curve $F_{\text{crit}}$ in the space $(\Omega,F)$.

For the  maximum mechanism, equation
\begin{equation}
C(\theta^*, E_{\text{max}}) = 0
\end{equation} is written as follows
\begin{equation}
E_{\text{max}} + F G(E_{\text{max}}) - \Omega J(E_{\text{max}}) = 0,
\end{equation} by solving which one can obtain
\begin{equation}\label{eq:max:case:one}
F_{\text{crit}} =  \frac{J(E_{max})}{G(E_{\text{max}})}\Omega - \frac{E_{\text{max}}}{G(E_{\text{max}})}.
\end{equation}

Now, we proceed to numerical verifications.
\subsubsection*{Example}
\begin{figure}[H]
\centering
\includegraphics*[width=0.5\textwidth]{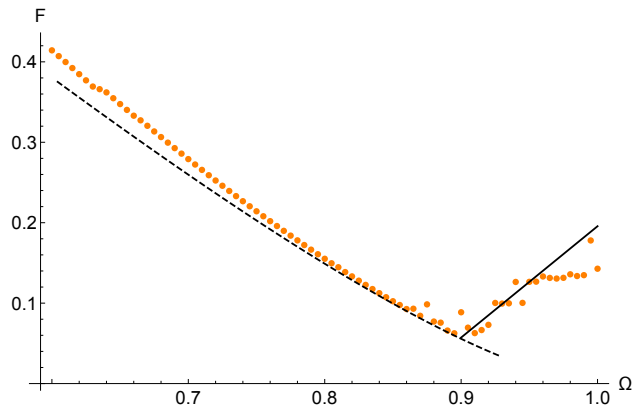}
\caption{Comparison of theoretical prediction of $F_{\text{crit}}(\Omega)$ (black lines) with numerical simulation (orange dots). The parameters are $\alpha = -1/2$ and $\beta = 1/20$. \label{fig:f:omega:plot}}
\end{figure}

\begin{figure}[H]
\begin{subfigure}[t]{0.32\textwidth}
\includegraphics*[width=0.9\textwidth]{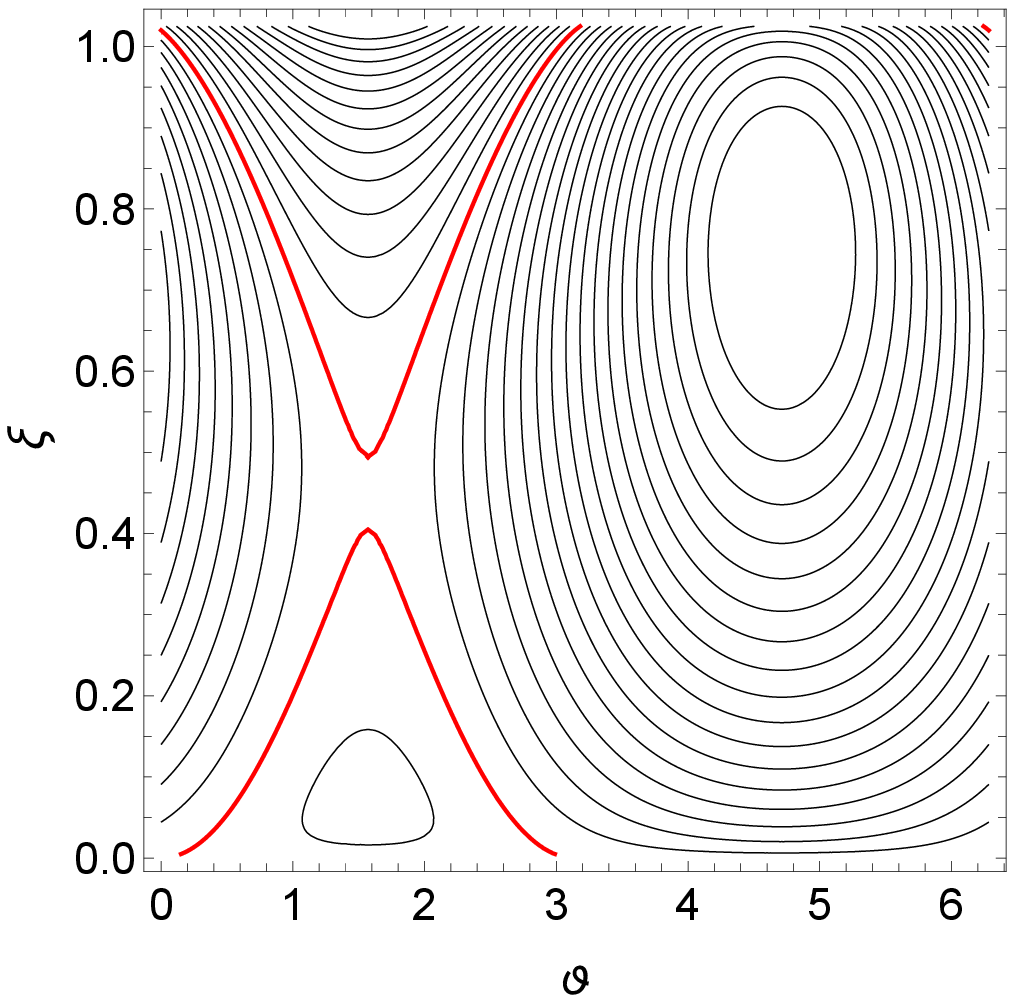}
\caption{}
\end{subfigure}
\hfill
\begin{subfigure}[t]{0.32\textwidth}
\includegraphics*[width=0.9\textwidth]{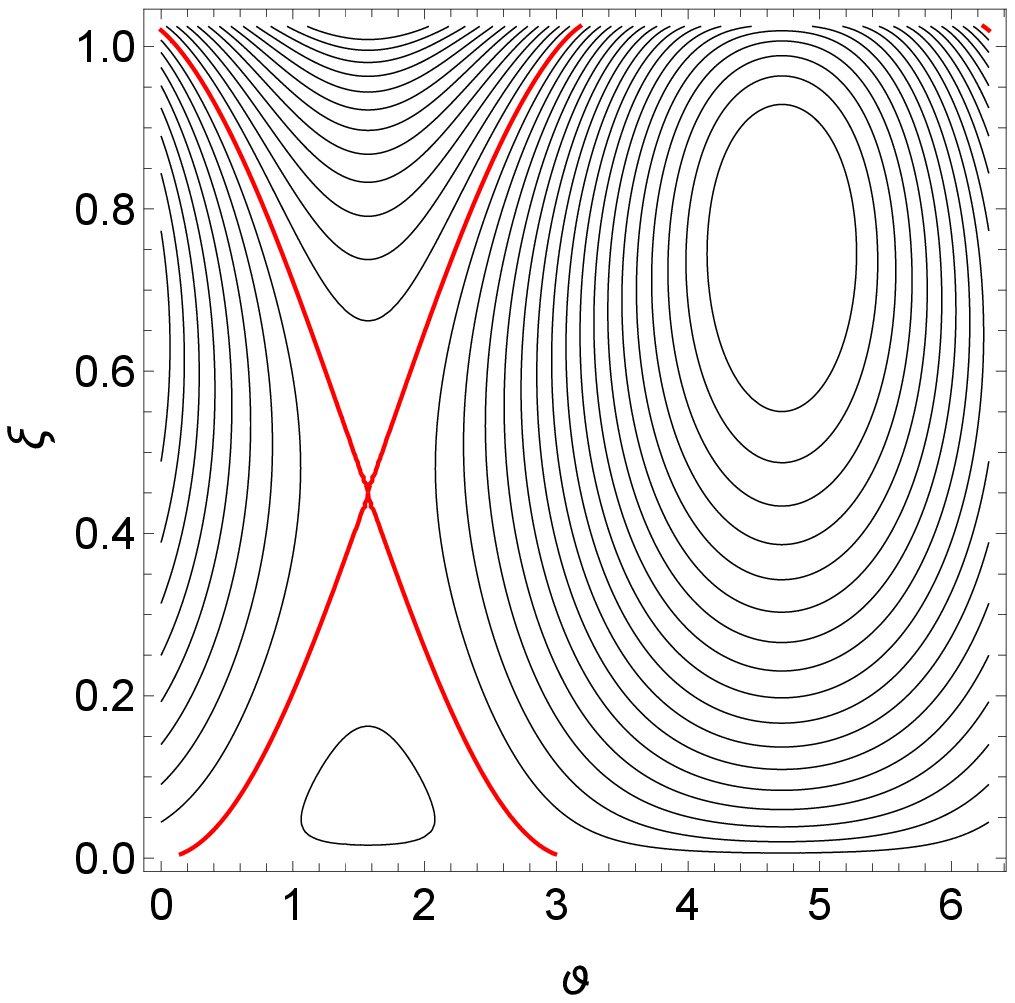}
\caption{}
\end{subfigure}
\hfill
\begin{subfigure}[t]{0.32\textwidth}
\includegraphics*[width=0.9\textwidth]{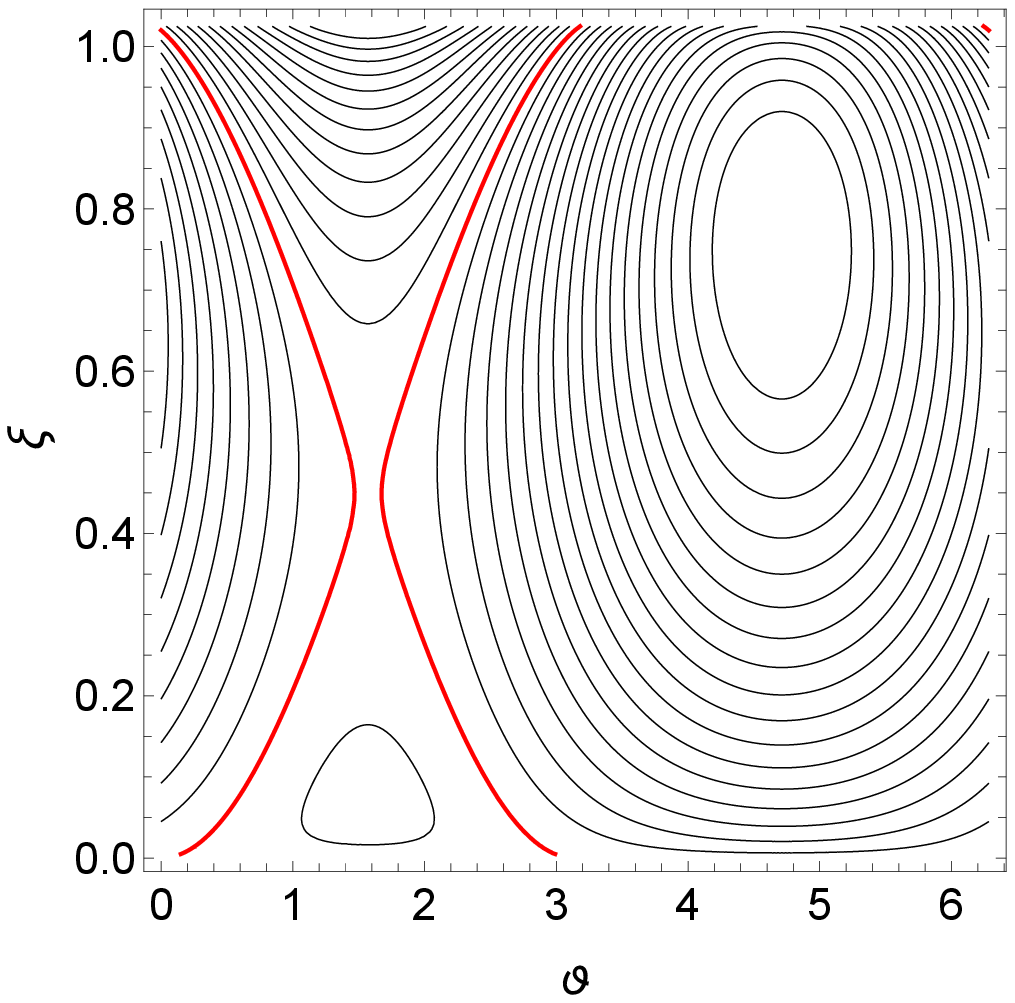}
\caption{}
\end{subfigure}
\caption{\label{fig:saddle:1}Transformation of the LPT (thick red curve) via the saddle mechanism as the amplitude $F$ of the external forcing passes the critical value~$F_{\text{crit}}$. Panels (a), (b) and (c) correspond to $F=0.099$, $0.09946$ and $0.1$, respectively. The values of other parameters are the following: $\Omega = 0.85$, $\alpha = -1/2$ and $\beta = 1/20$.}
\end{figure}

\begin{figure}[H]
\centering
\begin{subfigure}[t]{0.48\textwidth}
\centering
\includegraphics*[width=0.6\textwidth]{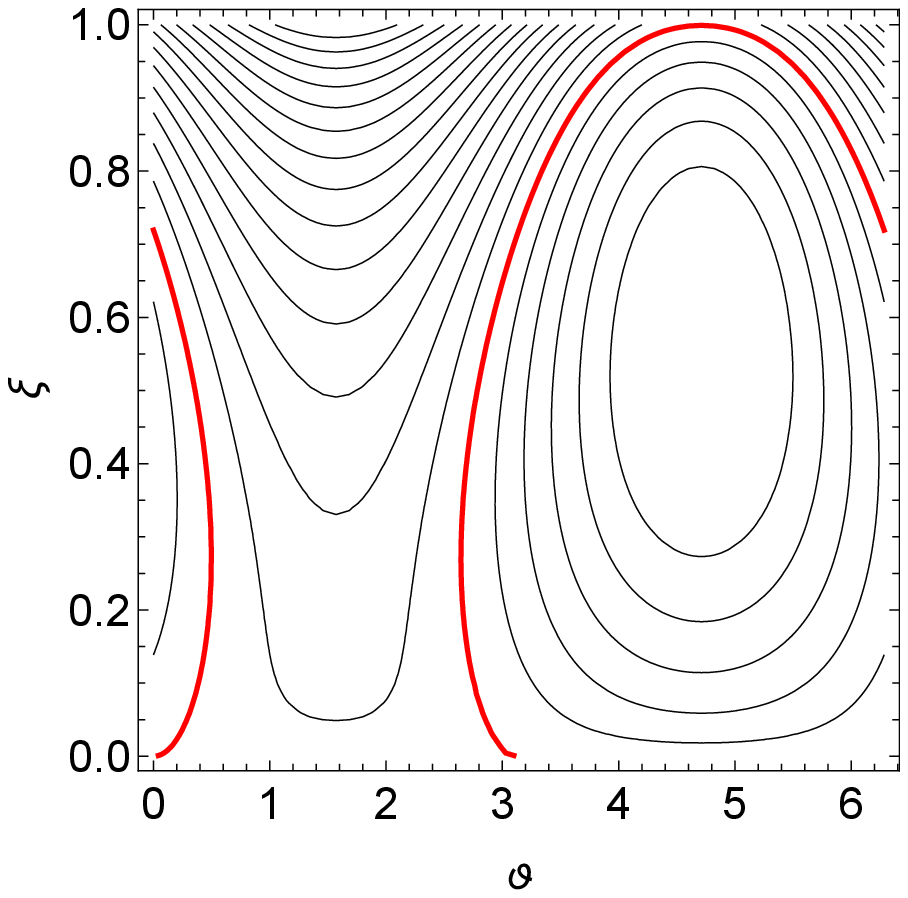}
\caption{}
\end{subfigure}
\begin{subfigure}[t]{0.48\textwidth}
\centering
\includegraphics*[width=0.6\textwidth]{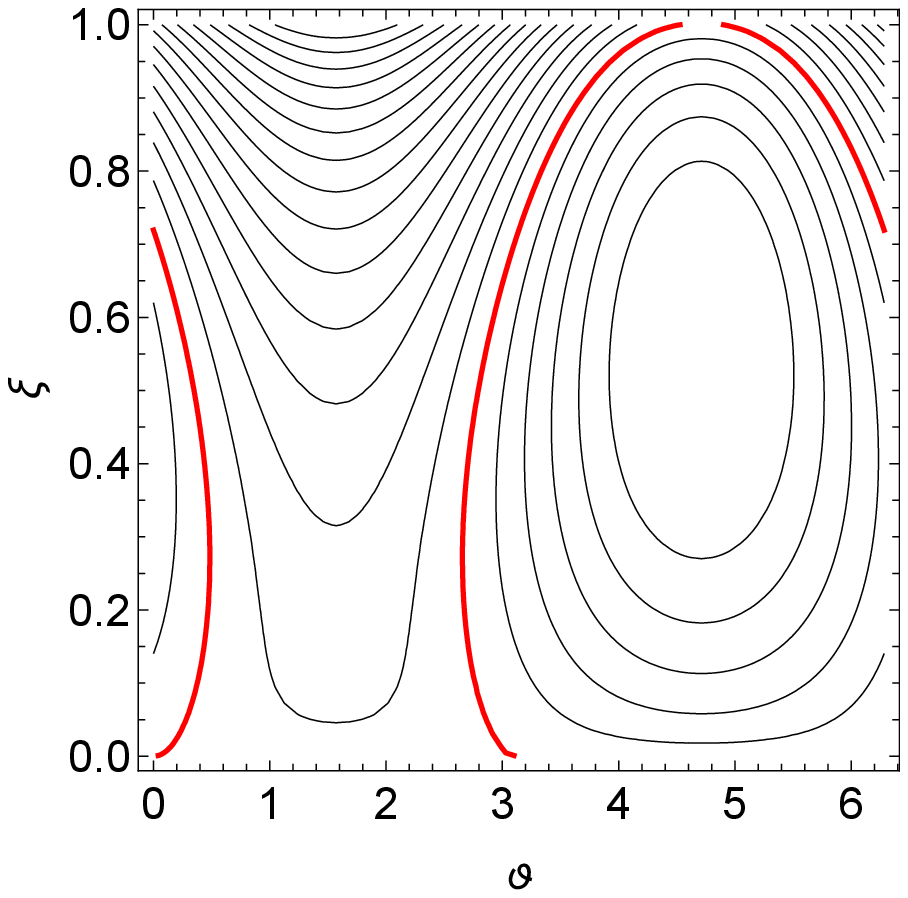}
\caption{}
\end{subfigure}
\caption{\label{fig:max:1}Illustration of the maximum mechanism for the excitation frequency $\Omega=0.92$. Pavel (a): $F =  0.0845$; panel (b): $F =  0.0865$. The notation and the other parameters are the same as in Figure~\ref{fig:saddle:1}.}
\end{figure}

\noindent
To illustrate formulae~\eqref{eq:saddle:case:one} and \eqref{eq:max:case:one} we select $\alpha = -1/2$ and $\beta = 1/20$ (see Figure~\ref{fig:case1:pot}). In this case $q_{\text{thres}}=2.76393$ and $E_{\text{max}} = 1.03006$.
Figure~\ref{fig:f:omega:plot} shows $F_{\text{crit}}(\Omega)$ near $1:1$ resonance. Dashed curve is obtained through the saddle mechanism, i.e., it is the graph of the parametric curve~\eqref{eq:saddle:case:one}. Solid line is defined by the maximum mechanism, i.e., it is a linear function $F_{\text{crit}}(\Omega)$ defined by equation~\eqref{eq:max:case:one}. Orange dots represent the results of numerical simulations.
Figures~\ref{fig:saddle:1},~\ref{fig:max:1} show level curves of the conservation law~\eqref{eq:cons:law:xi} and illustrate the saddle and the maximum mechanisms, respectively. Three panels of Figure~\ref{fig:saddle:1} portray the transformation of the phase cylinder as the amplitude~$F$ of the external forcing crosses a critical value $F \approx 0.0995$. Red curve represents the LPT. Equation $\frac{\partial C}{\partial \vartheta} = 0$  yields a saddle point $\left(\vartheta^\dagger, \xi^\dagger\right)$ where $\vartheta^\dagger=\pi/2$ and $\xi^\dagger$ can be used implicitly to parameterize the curve $F_{\text{crit}}(\Omega)$.
Similarly, the maximum mechanism is demonstrated on the Figure~\ref{fig:max:1}. Here LPT becomes tangent to the circle $\xi = E_{\text{max}}$ at $\vartheta^\dagger=3\pi/2$ when $F\approx 0.085$.

\begin{figure}[H]
\centering
\begin{subfigure}[t]{0.48\textwidth}
\centering
\includegraphics*[width=0.9\textwidth]{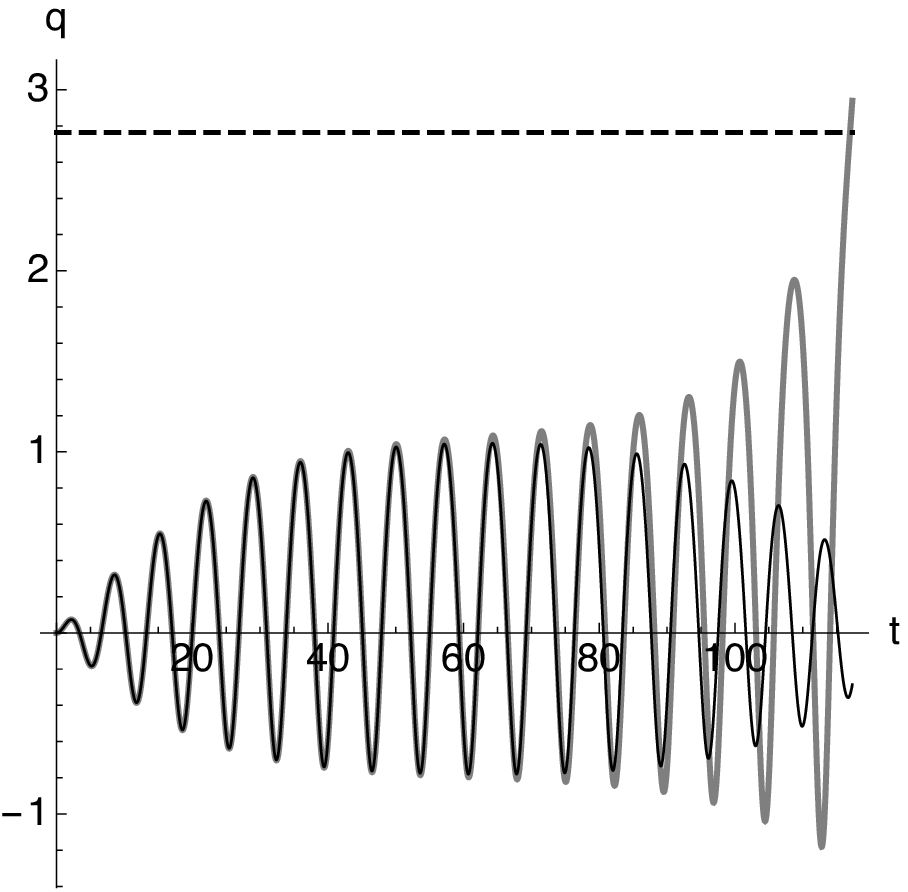}
\caption{}
\end{subfigure}
\hfill
\begin{subfigure}[t]{0.48\textwidth}
\centering
\includegraphics*[width=0.9\textwidth]{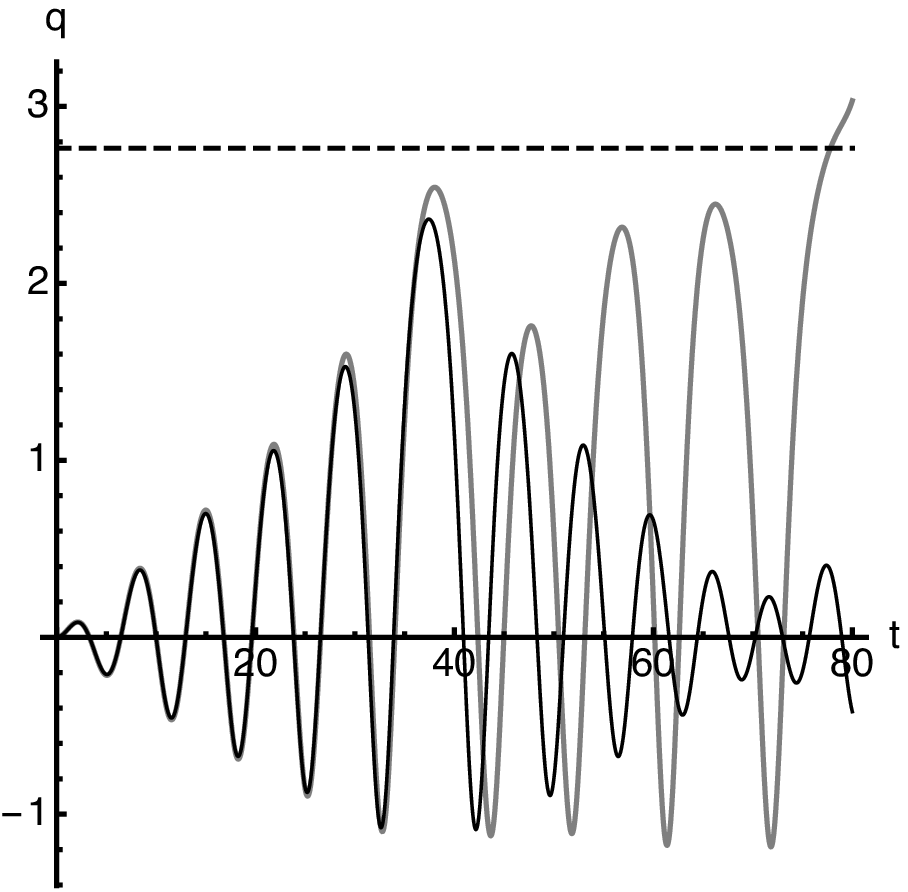}
\caption{}
\end{subfigure}
\caption{Time series of trajectories with initial conditions $q(0)=0$, $p(0)=0$ for different values of the parameter~$F$. Dashed line corresponds to the threshold value $q_{\text{thres}}=2.7639$. Panel (a): $F=0.0709$ (black), $F=0.0711$ (gray), $\Omega=0.88$; panel (b): $F=0.0845$ (black), $F=0.0846$ (gray), $\Omega=0.92$.\label{fig:timetraces}}
\end{figure}

Figure~\eqref{fig:timetraces} illustrates the difference between two mechanisms in terms of the time traces. The left panel shows time traces of two trajectories of system~\eqref{eq:main} starting from the zero initial condition. One of them undergoes escape however the other stays safely inside the potential well. Both averaged energy and the amplitude of oscillations undergo a drastic change. However, it's different from the maximum mechanism (see panel (b)). Here, both trajectories with the parameter~$F$ below and above the critical value have a commensurate value of averaged energy.

\subsection{Case II: Inverted Quartic Potential}
\noindent As we mentioned before there is no conceptual difference between two cases. For the sake of completeness, we briefly present main results and illustrate them with an example.
When $\beta$ is negative the potential has two maxima at 
\[
q = \frac{-\alpha\pm\sqrt{\alpha ^2-4 \beta } }{2 \beta }
\]
and a minimum at $q=0$ (bottom of the well), and $V(q)\to - \infty$ as $q \to \pm\infty$. The barrier $E_{\text{max}}$ is the smaller of two maxima. One can easily show that $E_{\text{max}}$ is
\[
E_{\text{max}} = \begin{cases}
-\frac{\left(\sqrt{\alpha ^2-4 \beta }+\alpha \right)^2 \left(\alpha  \left(\sqrt{\alpha
   ^2-4 \beta }+\alpha \right)-6 \beta \right)}{96 \beta ^3}, & \alpha > 0,\\
\frac{\left(\alpha -\sqrt{\alpha ^2-4 \beta }\right)^2 \left(\alpha  \left(\sqrt{\alpha
   ^2-4 \beta }-\alpha \right)+6 \beta \right)}{96 \beta ^3}, & \alpha < 0.
\end{cases}
\]

\noindent By following the same steps and notation as in as in Subsection 3.1, we arrive at the conservation law~\eqref{eq:cons:law:xi} for the case of inverted quartic potential:
\begin{equation}
C(\vartheta, \xi) = \xi -\frac{F\pi\sqrt{(b-d)(c-a)} \sinh(2\omega)}{2\bfK(k)\sinh(2\omega_0)}\sin{\vartheta} - \Omega J = \text{const},
\end{equation} where there averaged action $J$ is
\begin{align*}
&J=\frac{\sqrt{|\beta|}}{24 \pi  \sqrt{(a-c) (b-d)}} \left(\bfE(k) (a-c) (b-d) \left(3 a^2-2 a (b+c+d)+3 b^2-\right.\right.\\&\left.\left. 2 b
   (c+d)+3 c^2-2 c d+3 d^2\right)- (c-d) \left(\bfK(k) (d-b) \left(a^2+a (-4 b-4 c+6 d)\right.\right.\right.\\&\left.\left.\left.+3 b^2-2 b c+3 c^2-3 d^2\right)+3 \bfPi \left(\gamma ^2, k\right) \left(a^3-a^2
   (b+c+d)-\right.\right.\right.\\&\left.\left.\left.a \left(b^2-2 b (c+d)+(c-d)^2\right)+b^3-b^2 (c+d)-b (c-d)^2+(c-d)^2
   (c+d)\right)\right)\right),
\end{align*} where
\[
k = \sqrt{\frac{(b-c)(a-d)}{(a-c)(b-d)}}, \quad \gamma^2 = \frac{b-c}{b-d}.
\] Again, from this point we proceed to numerical simulations.
\subsubsection*{Example}
\begin{figure}[H]
\centering
\includegraphics*[width=0.5\textwidth]{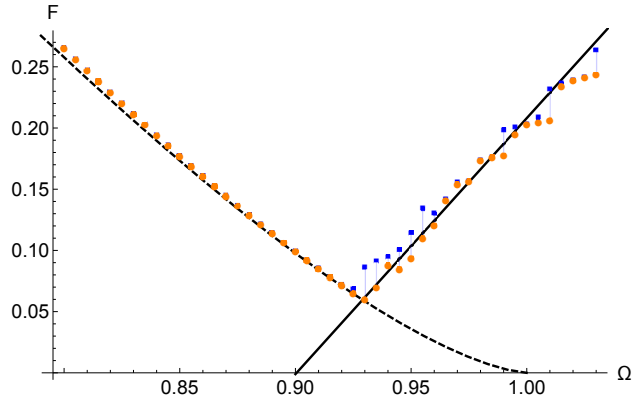}
\caption{Comparison of theoretical prediction of $F_{\text{crit}}(\Omega)$ with numerical simulations. Blue square markers represent simulations with the time limit of 100 periods; orange circles represent simulations with time limit 1000 periods. Dashed and solid curves are obtained theoretically through the saddle and the maximum mechanisms, respectively. The values of the parameters are $\alpha=-3/50$, $\beta=17/250$. \label{fig:f:omega:plot:2}}
\end{figure}

In order to illustrate Case II, we select the parameters to be $\alpha=-3/50$, $\beta=17/250$. Similarly, we obtain approximation for the curve $f_{\text{crit}}(\Omega)$ (see Figure~\ref{fig:f:omega:plot:2}).

\begin{figure}[H]
\begin{subfigure}[t]{0.32\textwidth}
\includegraphics*[width=\textwidth]{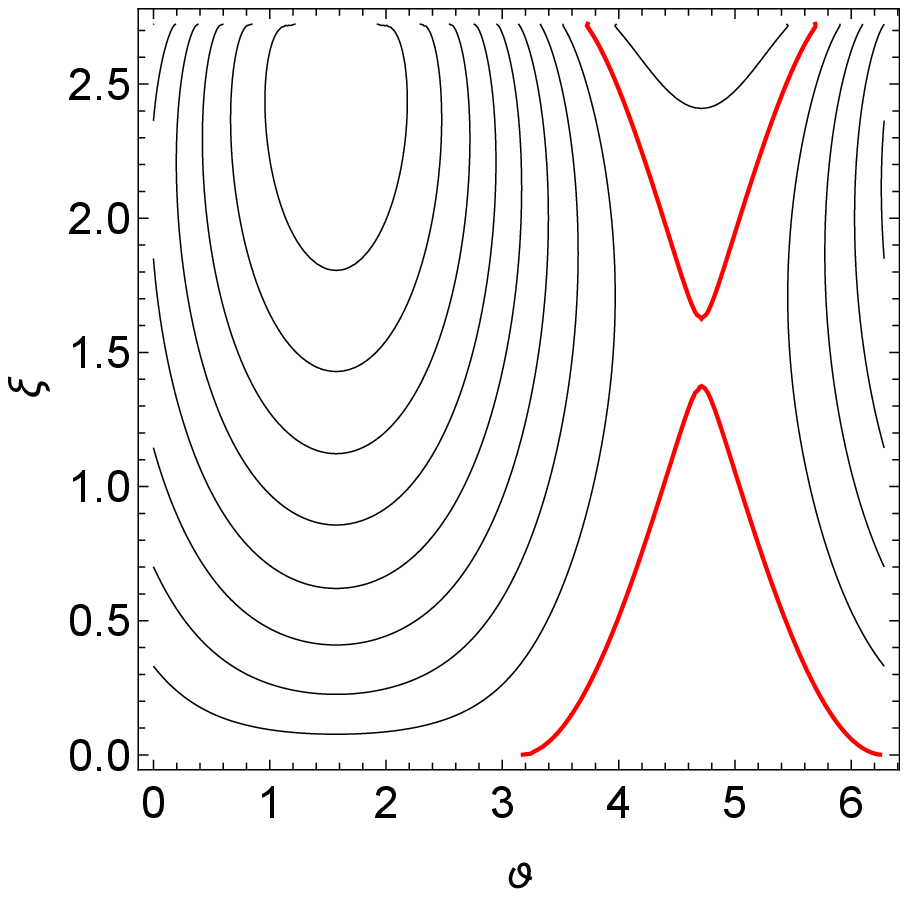}
\caption{}
\end{subfigure}
\hfill
\begin{subfigure}[t]{0.32\textwidth}
\includegraphics*[width=\textwidth]{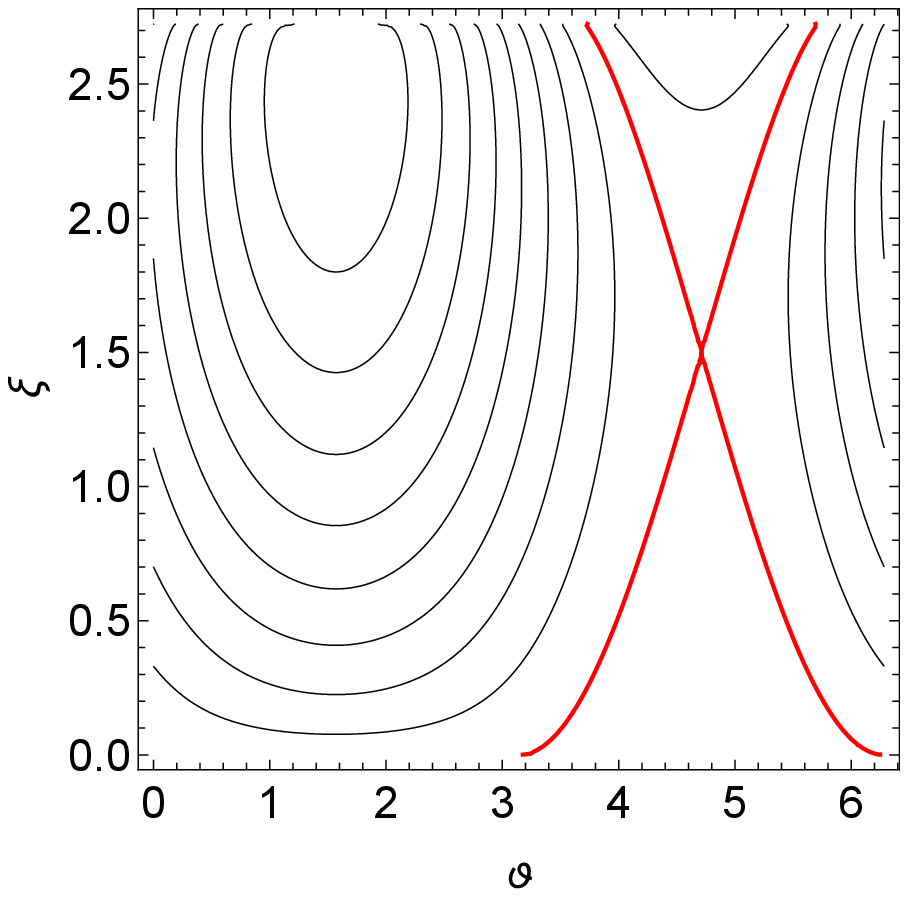}
\caption{}
\end{subfigure}
\hfill
\begin{subfigure}[t]{0.32\textwidth}
\includegraphics*[width=\textwidth]{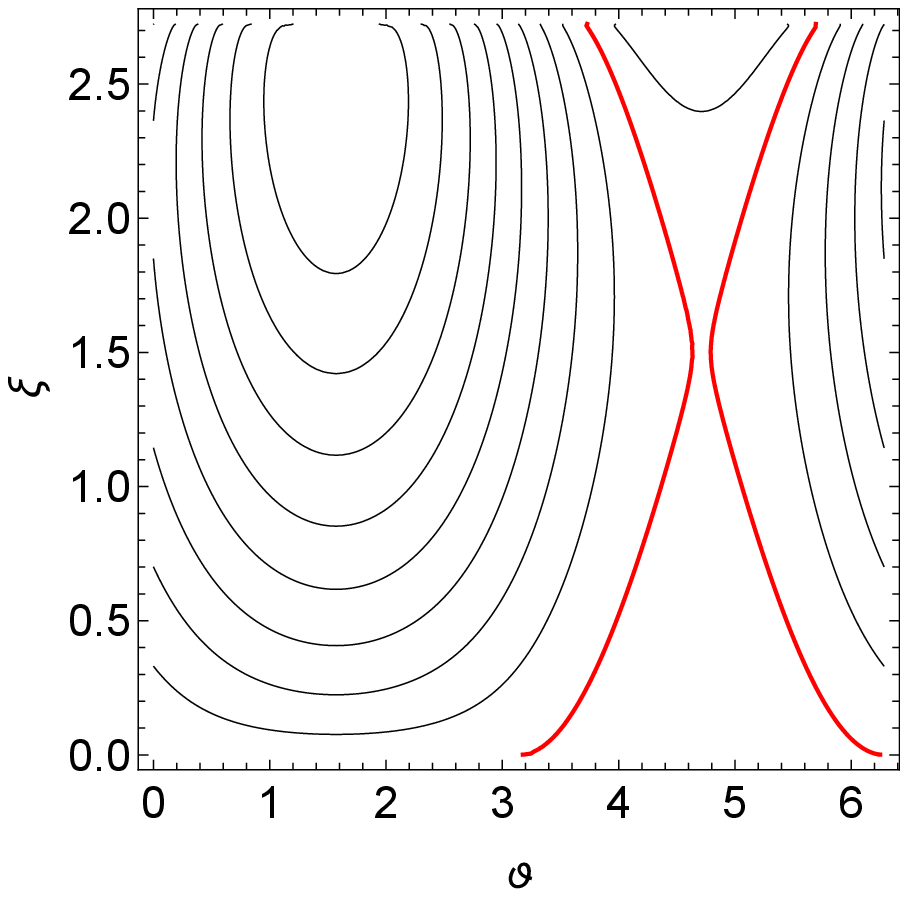}
\caption{}
\end{subfigure}
\caption{Saddle mechanism transition of the LPT (thick red curve) for $\Omega = 0.81859$, $\alpha=-3/50$, $\beta=17/250$. Panels (a), (b) and (c) correspond to $F = 0.208$, $0.20873$ and $0.2094$, respectively.}
\end{figure}

\section{Approximation of the model electrostatic potential}\label{sec:mems}

\subsection{Description of the Model}

\begin{figure}[H]
\centering
\includegraphics*[width=0.5\textwidth]{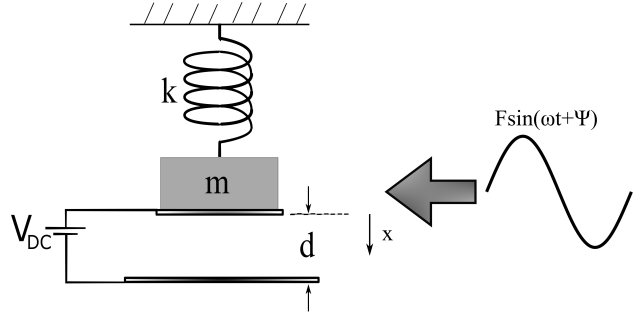}
\caption{\label{fig:mems}A sketch of a parallel-plate electrostatic actuator under the influence of external harmonic force.}
\end{figure}

\noindent A simple MEMS device is a parallel-plate electrostatic actuator~\cite{younis2011mems}, see Figure~\ref{fig:mems}. The device consists of a parallel-plate capacitor with a moving upper electrode attached to a spring and a stationary lower electrode. A small DC load applied to the capacitor creates electrostatic force compensated by mechanical restoring force of the spring, thus, pulling the plate in a new equilibrium position. One can increase the voltage up to some critical value for which the restoring force of the spring cannot longer resist the opposing electrostatic force. Inevitably, this results in the plates collapsing. The described phenomenon is a structural instability called {\it static pull-in}. In applications to resonators, AC load is applied in addition to the DC load in which case the pull-in can occur at much smaller values of the critical DC voltage. If a pull-in happens due to the AC loading, it is called a {\it dynamic pull-in}. Alternatively, one can consider a plate excited by an external mechanical vibration.

Regardless, the equation of motion of the upper plate of mass $m$ without a damping under the influence of external harmonic force is 
\begin{equation}
\label{eq:mems:dim}
m\ddot{x}+k x = \frac{\varepsilon A V^2_{\text{DC}}}{2(d-x)^2} + f \sin{(\omega t+\Psi)},
\end{equation} where $k$ is the spring coefficient, $d$ is the gap width, $\epsilon$ is the dielectric constant of the gap medium, $V_{\text{DC}}$ is the input voltage, $A$ is the area of the electrode. The external harmonic force can be due to additional AC loading, or to the external mechanical vibrations.

\begin{figure}[H]
\centering
\includegraphics*[width=0.6\columnwidth]{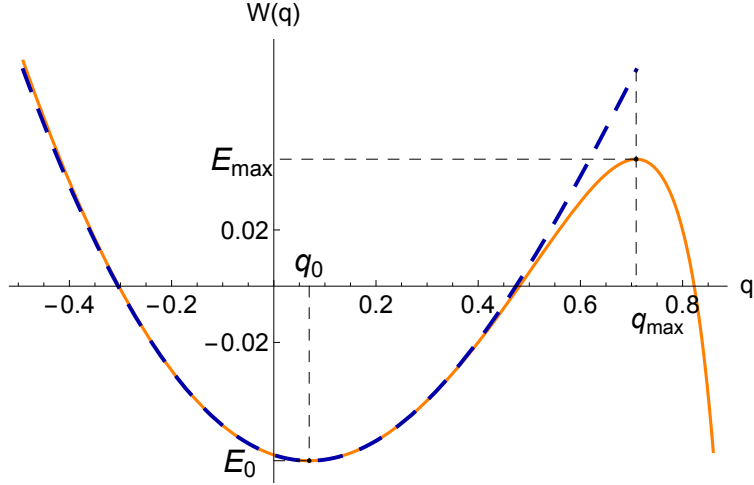}
\caption{Orange solid curve is a graph of potential $W\left(q\right)$ with $\nu=0.06$ and $d=1$. Blue dashed curve is a graph of the approximating potential~$\widetilde{W}$.\label{fig:mems-potentials}}
\end{figure}

By introducing rescaled time $\tau = \sqrt{\frac{k}{m}} t$ and putting $q(\tau) = x(t)$, one can rewrite equation~\eqref{eq:mems:dim} as follows
\begin{equation}
q''+q=\frac{\nu}{\left(d-q\right)^2}+F\sin{\left(\Omega\tau+\Psi\right)},
\end{equation} where
\[
\nu=\frac{\epsilon A\ V_{DC}^2}{2m}, \quad F=\frac{f}{m} \quad \mathrm{and}\quad\Omega=\omega\sqrt{\frac{m}{k}},
\]
and, the prime symbol ($'$) denotes differentiation with respect to the rescaled time $\tau$.
The potential energy of the unforced system is
\begin{equation}\label{eq:mems:potential}
W\left(q\right)=\frac{q^2}{2}-\frac{\nu}{d-q}.
\end{equation}
An example of potential $W\left(q\right)$ is presented by a blue curve on Figure~\ref{fig:mems-potentials}. The escape occurs when $q\left(t\right)$ crosses the threshold value $q=q_{\mathrm{max}}$. The initial conditions $q(0)=q_{0}$, $q^{\prime}(0)=0$ correspond to the minimum of energy $E=E_0$.

We want to find the minimal amplitude $f_{\mathrm{crit}}\left(\Omega\right)$ of external forcing needed for the escape. For a given potential~\eqref{eq:mems:potential} it is a difficult if not impossible problem. That is why we suggest several forth order polynomials as candidates for approximation. The approaches we take can be classified into two types: global approximation and local approximation. The global approximation is an ad hoc approach to fit a polynomial curve onto a given potential. The local approximation utilizes the Taylor's polynomial near the minimum.

\subsection{Global approximation}

\noindent The idea behind the global approximation is to approximate the given electrostatic potential~$W$ with a handful parameters such as its height, width and the curvature at the minimum. Alas, it does not work. In fact, it is difficult if not impossible to find a valid approximation using so little data. In addition to the three parameters listed above, one has to take into consideration the other side of the well, curvature at the maximum, etc.

Before approximating potential~$W$ it is useful to introduce a translated potential~$\widehat{W}$ with the minimum exactly at the origin:
\[
\widehat{W}(q) = W(q+q_0)-W(q_0).
\]

Then, the modified threshold energy level becomes $\widehat{E}_0=E_0-W(q_0)$. The motion of the particle inside the potential~$\widehat{W}$ is analogous to dynamics inside $W$ modulo the coordinate translation $q \mapsto q - q_0$.

\begin{figure}[H]
\centering
\includegraphics*[width=0.6\columnwidth]{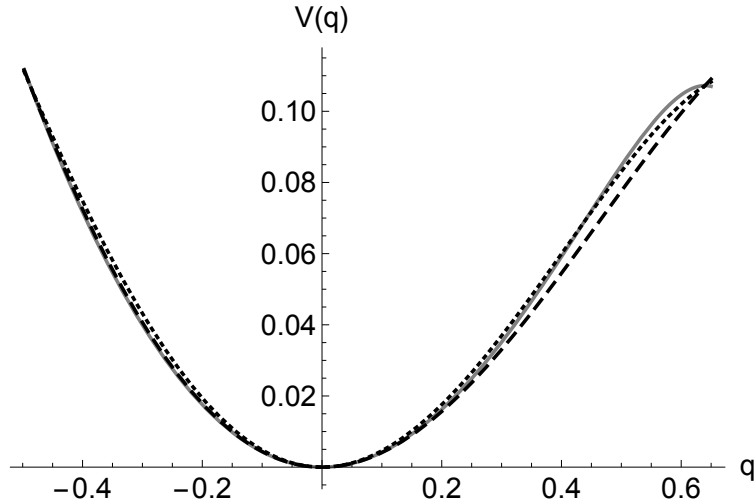}
\caption{Gray, dashed and dotted curves correspond to the graphs of $\widehat{W}$, $p_1$ and $p_2$, respectively.\label{fig:mems:comparison}}
\end{figure}

As an example we consider two approximations. Let $p$ be a forth-order polynomial function:
\[
p(x) = \frac{a}{2} x^2 + \frac{b}{3} x^3 + \frac{c}{4} x^4.
\]
The first approximation (orange on Figure~\ref{fig:mems:comparison}) is obtained by solving the following equations
\begin{equation}\label{eq:global:appr:1}
p(q_{\text{max}})=\widehat{W}(q_{\text{max}}),\quad p(q_{\text{min}})=\widehat{W}(q_{\text{min}}),\quad p''(0)=\widehat{W}''(0).
\end{equation}
Similarly, for the second approximation (green) the coefficients $a$, $b$, $c$ are the solution to
\begin{equation}\label{eq:global:appr:2}
p(q_{\text{infl}})=\widehat{W}(q_{\text{infl}}),\quad p(q_{\text{min}})=\widehat{W}(q_{\text{min}}),\quad p''(0)=\widehat{W}''(0)
\end{equation} where $q_{\text{infl}}$ corresponds to the inflection point of $\widehat{W}$, i.e., solution to $\widehat{W}''(q_{\text{infl}})=0$.

For example, if $\nu=0.06$ and $\delta=1$, then the boundaries of the well are $q_{\text{max}}=0.639856$ and $q_{\text{min}}=-0.487499$. By solving equations~\eqref{eq:global:appr:1} one obtains a polynomial:
\[
p_1(x) = 0.475492 x^2 -0.116032 x^3 - 0.340386 x^4.
\] Likewise, equations~\eqref{eq:global:appr:2} yields the second global approximation polynomial:
\[
p_2(x)=0.425583 x^2 - 0.140409 x^3 - 0.180387 x^4.
\]

\begin{figure}[H]
\centering
\includegraphics*[width=0.6\columnwidth]{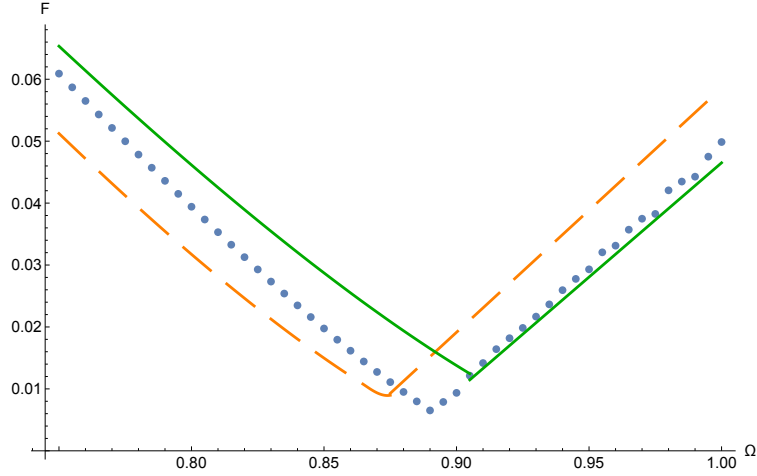}
\caption{Comparison of the critical forcing curves as functions of the frequency~$\Omega$. Blue dots represent the numerically obtain critical force curve for the potential $\widehat{W}$; dashed orange and solid green curves correspond to the escape from the approximating potentials $p_1$ and $p_2$, respectively.\label{fig:mems:global:appr:escape}}
\end{figure}

Both functions $p_1$ and $p_2$ superimposed onto the potential~$\widehat{W}$ are presented on Figure~\ref{fig:mems:comparison}. The corresponding critical escape curves in the parameter space $(\Omega, F)$ are shown on Figure~\ref{fig:mems:global:appr:escape}.

As one can see, the results are very sensitive to the initial form of the potential we choose. Two visually same approximations yield substantially different critical escape curves in the parameter space $(\Omega, f)$.

\subsection{$L^2$-heuristic approach}

\noindent Another way to obtain an approximation is to seek a truncated forth-order polynomial $p(x)$ that minimizes the following functional:
\begin{equation}
\int\limits_{q_{\text{min}}}^{q_{\text{max}}}\left[\widehat{W}(x)-p(x)\right]^2 dx.
\end{equation}

The minimizing polynomial $p(x)$ follows function $V(x)$ on the interval $\left[q_{\text{min}},\;q_{\text{max}}\right]$ and therefore, it is a good candidate for an approximating potential for the escape problem.

\begin{figure}[H]
\centering
\includegraphics*[width=0.5\columnwidth]{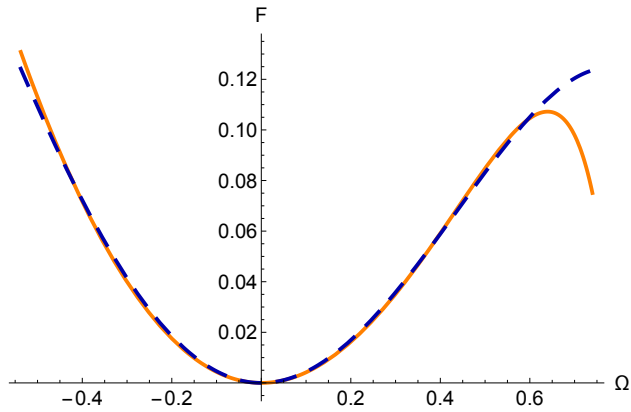}
\caption{Orange solid curve is potential $\widehat{W}(x)$ and blue dashed curve is approximating potential $p_3(x)$\label{fig:l2:approx}}
\end{figure}

Again with the chosen parameters $\nu=0.06$ and $\delta = 1$, the minimizing polynomial becomes
\[
p_3(x)=0.45386 x^2-0.103971 x^3-0.276211 x^4.
\]
Both functions $\widehat{W}(x)$ and $p_3(x)$ are plotted on Figure~\ref{fig:l2:approx}.

\begin{figure}[H]
\centering
\includegraphics*[width=0.5\columnwidth]{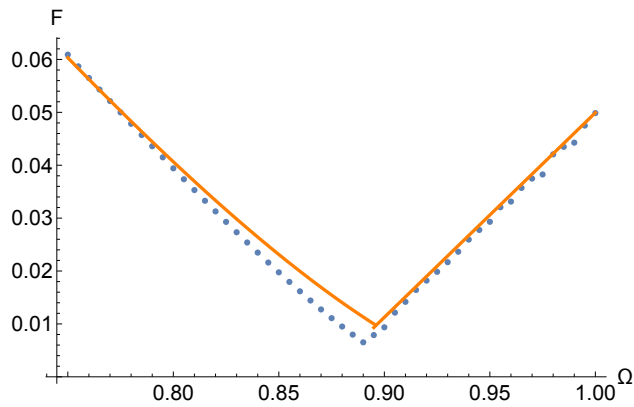}
\caption{Critical forcing amplitude as a function of the frequency $\Omega$.\label{fig:l2:escape}. Blue dots correspond to numerical values obtained for the potential~$\widehat{W}$. Orange curves represent the analytic prediction of the escape curve for the approximating polynomial~$p_3$.}
\end{figure}

The corresponding $F_{\text{crit}}$ curves are depicted on Figure~\ref{fig:l2:escape}.

\subsection{Local approximation}

The local approximation is the following function:
\begin{equation}\label{eq:local:approx}
\widetilde{W}\left(q\right)=E_0+\alpha_1\left(q-q_0\right)^2+\alpha_2\left(q-q_0\right)^3+\alpha_3\left(q-q_0\right)^4,
\end{equation} where
\begin{equation*}
\alpha_1=\frac{1}{2}-q_0^{\frac{3}{2}}\nu^{-\frac{1}{2}},\quad\alpha_2=-q_0^2\nu^{-1},\quad\alpha_3=-q_0^{\frac{5}{2}}\nu^{-\frac{3}{2}},
\end{equation*}
truncated at the energy level $E=\widetilde{E}_{\mathrm{max}}:=\widetilde{W}(q_{\text{max}})$. In other words, we approximate the potential energy~\eqref{eq:mems:potential} by taking its Taylor's expansion near $q=q_0$ up to the forth order term.  Note that $q_0=\nu/d^2+O\left(\nu^2\right)$, therefore, for small values of $\nu$, function~$W\left(q\right)$ is a weakly nonlinear potential well, i.e., $\alpha_2, \alpha_3 \ll\alpha_1$.

\begin{figure}[H]
\centering
\includegraphics*[width=0.5\textwidth]{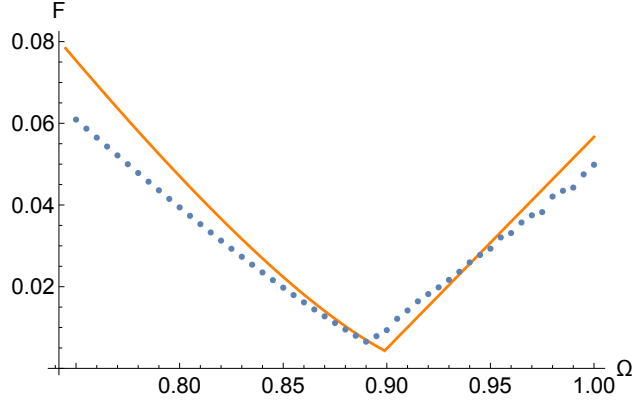}
\caption{Comparison of numerically obtained values of critical forcing $F_{\text{crit}}$ for $W(q)$ (blue dots), as well as theoretical prediction for the approximating polynomial $\widetilde{W}(q)$ (orange curve).\label{fig:mems:escape}}
\end{figure}

The comparison of theoretic prediction of $f_{\text{crit}}$ for the approximating quartic potential and numerical simulations for the electrostatic potential is presented on Figure~\ref{fig:mems:escape}. As we can see, the proposed method yields a decent approximation of the~$f_{\text{crit}}(\Omega)$ for the exact potential.

One can observe that despite the fact that most of the discrepancy between the potential and its approximation occurs near the right edge of the well, it significantly impacts the escape curve. In particular, it effects the position of the minimum corresponding to the resonance frequency.

\subsubsection{Comparison of approximation orders}

\begin{figure}[H]
\centering
\begin{subfigure}[t]{0.32\textwidth}
\includegraphics*[width=\columnwidth]{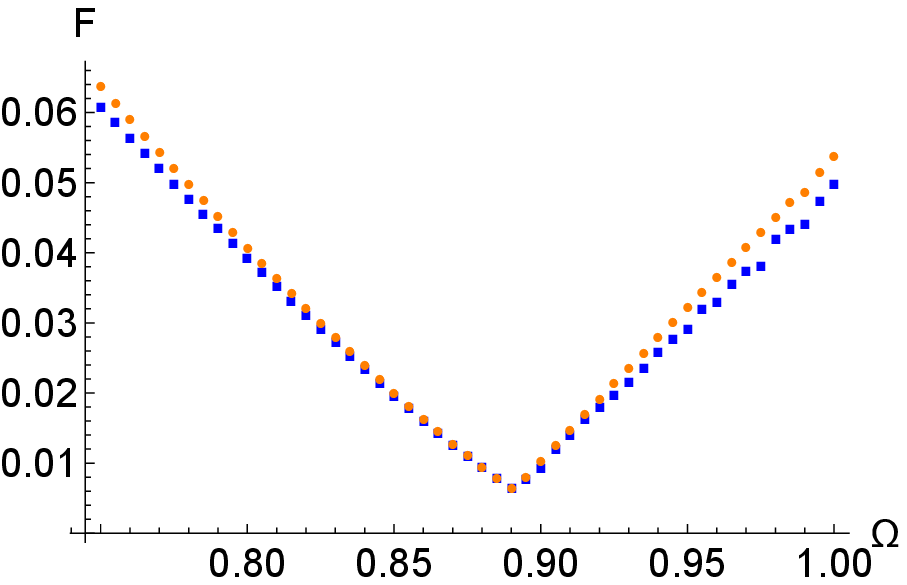}
\caption{}
\end{subfigure}
\hfill
\begin{subfigure}[t]{0.32\textwidth}
\includegraphics*[width=\columnwidth]{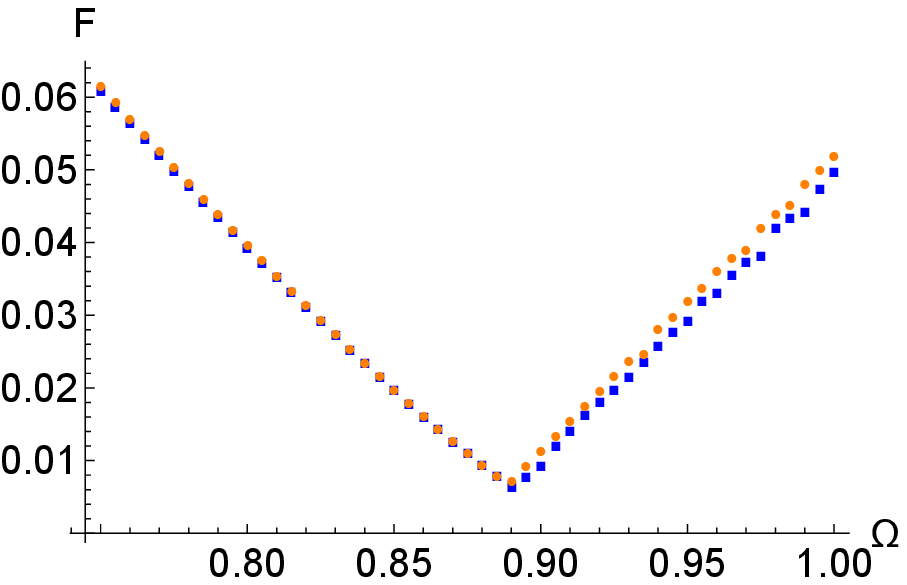}
\caption{}
\end{subfigure}
\hfill
\begin{subfigure}[t]{0.32\textwidth}
\includegraphics*[width=\columnwidth]{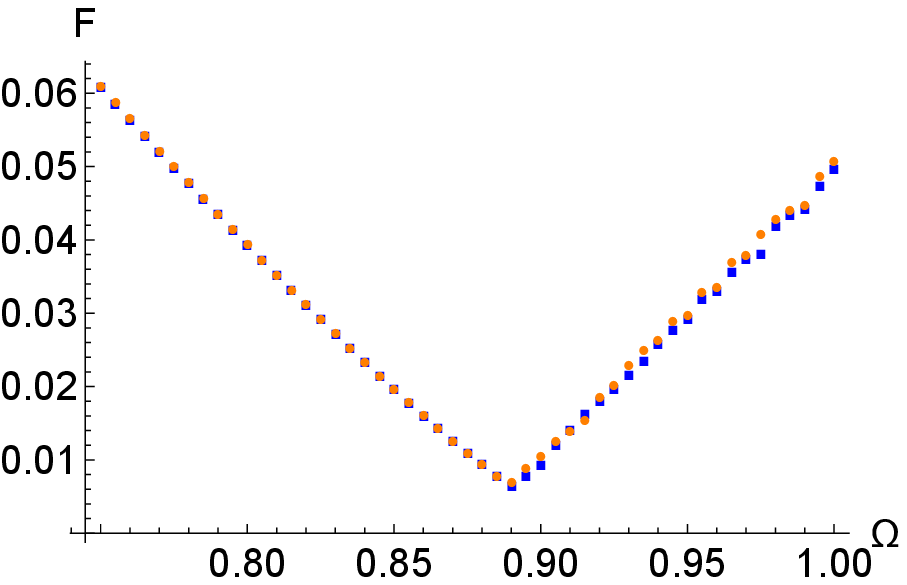}
\caption{}
\end{subfigure}
\caption{Comparison of $F(\Omega)$ higher-order local approximations. Blue square markers correspond to the critical escape values for the potential~\eqref{eq:mems:potential}, orange circles denote the approximation. Panels (a), (b), (c) correspond to the 6th, 8th, and 10th order, respectively. \label{fig:approx:order}}
\end{figure}

\noindent In order to obtain a better approximation of the escape curve $F(\Omega)$ one can expand local approximation~\eqref{eq:local:approx} with higher order terms. Unfortunately, the analytical method presented in Section~\ref{sec:model:quartic} becomes inapplicable, as there is no AA representation for the polynomial potentials of order higher than four, Therefore, the further comparison is performed numerically. Three panels of Figure~\ref{fig:approx:order} show comparison of local approximations of order 6,8 and 10, respectively.

As expected, the quality of the approximation improves as the order increases.

\section{Conclusions}

\noindent The results presented above demonstrate that in the problem of forced escape the idea of approximating the realistic potential functions by tractable low-order polynomials  is in principle viable, but somewhat tricky. From one side, the V-shaped dependence of the escape threshold on the excitation frequency reveals itself in all approximation methods, both local and global. Moreover, the sharp minimum at this curve is predicted by all approximations with relative accuracy of at least 7-10 percent. Such accuracy can be considered as satisfactory, since the inaccuracy of the model potential, and especially the errors related to reduction to the single-mode approximation, can introduce much more severe errors. In addition, the time series of the response reveal that for the exact model potential one encounters the well-known mechanisms of escape in the conditions of 1:1 resonance (maximum mechanism and saddle mechanism), despite the fact that the RM cannot be presented in analytically explicit form.

From the other side, it is somewhat surprising that minor variations of the approximating potential, almost invisible to the eye, lead to quite noticeable modifications of the escape threshold curve. It points on a considerable sensitivity of the escape threshold to the details of the model. In reality, it might mean that statistical approach will be inevitable to get reliable information on possible range of the escape thresholds.

Among the methods presented in this work, it is worth noting the $L_2$-heuristic approximation which yields the best estimate for the escape curve comparing to the other approaches. Local approximation based on the Taylor's polynomial near the minimum is another viable technique. The quality of local approximation increases with the order of the polynomial. Unfortunately, the analytic prediction cannot be obtained in the framework of the described general approach for any polynomial of order higher than four, at least in terms of elliptic functions.

\section*{Acknowledgements}

\noindent The authors are very grateful to Israel Science Foundation (grant 1696/17) for financial support.

\section*{Appendix}
In the Appendix we present the derivations of the transformation to AA variables and the conservation law~\eqref{eq:cons:law} of the slow-flow equations for the quartic potential. For the sake of brevity,  we restrict ourself only to Case I, i.e., $V(q)$ is a double-well potential. All the derivations for the inverted quartic potential (Case~II) are completely analogous.
According to~\eqref{eq:aa}  the action variable is
\begin{equation}\label{eq:I}
\begin{aligned}
I(E) = \frac{1}{2\pi}\oint\limits_{\Gamma_E} p(q,\, E)\mathrm{d}q &= \frac{\sqrt{2}}{\pi}\int\limits_{d}^{c} \sqrt{E-\frac{q^2}{2}-\alpha \frac{q^3}{3}-\beta\frac{q^4}{4}} \mathrm{d}q \\&= \frac{\sqrt{2\beta}}{2\pi}\int\limits_d^c\sqrt{(a-q)(b-q)(c-q)(q-d)}\mathrm{d}q
\end{aligned}
\end{equation}  where $a>b>c>d$ are the roots of the forth-order polynomial equation
\[
E-\frac{q^2}{2}-\alpha \frac{q^3}{3}-\beta\frac{q^4}{4} = 0.
\]

The last integral in~\eqref{eq:I} is a table integral (see~\cite{byrd2013handbook}) expressed as follows
\begin{align*}
I = &\frac{1}{48\pi}\sqrt{\frac{2\beta}{(a-c)(b-d)}}\left[(a-c)(b-d)\left(\frac{16 \left(\alpha ^2-3 \beta \right)}{3 \beta ^2}\right)\mathbf{E}(k) +\right.\\ 
&(a-c) (a-d) \left(3 a^2-6 a b-b^2+4 b (c+d)-3 c^2+2 c d-3 d^2\right)\mathbf{K}(k) + \\
&\left.3(a-d)\left(-3 a^3+\frac{16 \alpha ^2 a}{9 \beta ^2}-\frac{4 a (a \alpha +3)}{3 \beta}-b^3+(c+d) \left(b^2-(c-d)^2\right)+b \left(c^2+d^2\right)\right)\bfPi(\gamma^2, k)\right]
\end{align*} where
\[
k = \sqrt{\frac{(a-b) (c-d)}{(a-c) (b-d)}}, \qquad \gamma^2 = \frac{d-c}{a-c} < 0,
\] and $\mathbf{K}(k)$, $\mathbf{E}(k)$, $\bfPi(\gamma^2,\, k)$ are the complete elliptic integrals of the first, the second and the third kind, respectively.

The angle variable is
\[
\theta = \frac{\partial}{\partial I} \int\limits_d^q {p\left(x, I\right) \mathrm{d}x} = \Omega(I) \frac{\partial}{\partial E} \int\limits_d^q {p\left(x, E\right) \mathrm{d}x}, \qquad \Omega(I)=\frac{\mathrm{d}E}{\mathrm{d}I}.
\]
By the Inverse Function Theorem
\begin{equation}
\label{eq:recip:freq}
\begin{aligned}
\frac{1}{\Omega(E)} &= \frac{\mathrm{d}I}{\mathrm{d}E} = \frac{\sqrt{2}}{2\pi}\int\limits_d^c\frac{\mathrm{d}q}{\sqrt{E-V(q)}} \\ &= \frac{1}{\pi}\sqrt{\frac{2}{\beta}}\int\limits_d^c \frac{\mathrm{d}q}{\sqrt{(a-q)(b-q)(c-q)(q-d)}}\\& = \frac{1}{\pi}\sqrt{\frac{2}{\beta(a-c)(b-d)}}\,\mathbf{K}(k).
\end{aligned}
\end{equation}
Also,
\begin{equation}
\label{eq:angle:E}
\begin{aligned}
\frac{\partial}{\partial E} \int\limits_d^q {p\left(x, E\right) \mathrm{d}x} &= \frac{1}{\sqrt{2}}\int\limits_d^q {\frac{\mathrm{d}x}{\sqrt{E-V(x)}}} \\&= \sqrt{\frac{2}{\beta}}\int\limits_d^q \frac{\mathrm{d}x}{\sqrt{(a-x)(b-x)(c-x)(x-d)}} \\&= \sqrt{\frac{2}{\beta(a-c)(b-d)}}\, \mathbf{F}\left(\phi,\, k\right)
\end{aligned}
\end{equation} where
$\phi = \arcsin{\left(\sqrt{\frac{(a-c)(q-d)}{(c-d)(a-q)}}\right)}$ and $\mathbf{F}\left(\phi, \,k\right)$ is the incomplete elliptic integral of the first kind.

Combining~\eqref{eq:recip:freq} and \eqref{eq:angle:E}, one obtains an equation
\[
\theta = \frac{\pi\, \bfF\left(\phi,\, k\right)}{\bfK(k)}
\] solving which for $q$ results in the following expression
\begin{equation}
\begin{aligned}
q\left(\theta,\, E\right) &= \frac{d(a-c)+a(c-d)\,\sn^2{\left(\frac{\bfK(k)}{\pi}\theta,\, k\right)}}{a-c+(c-d)\,\sn^2{\left(\frac{\bfK(k)}{\pi}\theta,\, k\right)}}\\&
= a + \frac{d-a}{1-\gamma^2 \sn^2{\left(\frac{\bfK(k)}{\pi}\theta,\,k\right)}},
\end{aligned}
\end{equation}
where $\sn\left(\cdot,\cdot\right)$ is the Jacobi elliptic sine function.

Fourier expansion of $q(\theta,\, E)$ can be obtained using well-known formulae (see for example~\cite{langebartel1980fourier}):
\begin{equation}
\label{eq:q:fourier}
q(\theta,\, E) = a + \frac{(d-a)\,\bfPi\left(\gamma^2,\, k\right)}{\bfK\left(k\right)}-\frac{\pi  \sqrt{(b-d)(a-c)}}{\bfK\left(k\right)}\sum\limits_{n=1}^\infty{\frac{\sinh (2 n \omega)}{\sinh (2n\omega_0)}\cos(n\theta)}
\end{equation} where
\[
\omega = \frac{\pi \left(\bfK\left(k^\prime\right)-\nu\right)}{2 \bfK(k)},\qquad \omega_0 = \frac{\pi \bfK\left(k^\prime\right)}{2 \bfK(k)}, \qquad k^\prime = \sqrt{1-k^2}.
\] and $\nu$ is defined by
\[\text{cn}\left(\nu,\, k^\prime\right)=\sqrt{\frac{c-d}{a-d}}, \qquad 0 < \nu < \bfK(k^{\prime}).\]
In particular, coefficient $q_1$ is
\[
q_1 = \bar{q}_1 = -\frac{\pi\sqrt{(b-d)(a-c)} \sinh(2\omega)}{2\bfK(k)\sinh(2\omega_0)}.
\]
Therefore, the conservation law~\eqref{eq:cons:law:xi} becomes
\begin{equation}
C(\vartheta, \xi) = \xi -\frac{F\pi\sqrt{(b-d)(a-c)} \sinh(2\omega)}{2\bfK(k)\sinh(2\omega_0)}\sin{\vartheta} - \Omega J = C.
\end{equation}

\bibliography{biblio.bib}

\end{document}